\documentclass{article}%
\usepackage{amssymb}
\usepackage{amsmath}
\usepackage{amsfonts}
\usepackage{graphicx}%
\newtheorem{theorem}{Theorem}

\newtheorem{conjecture}[theorem]{Conjecture}
\newtheorem{corollary}[theorem]{Corollary}

\newtheorem{definition}[theorem]{Definition}
\newtheorem{example}[theorem]{Example}

\newtheorem{lemma}[theorem]{Lemma}

\newtheorem{proposition}[theorem]{Proposition}
\newtheorem{remark}[theorem]{Remark}

\newenvironment{proof}[1][Proof]{\noindent\textbf{#1.} }{\ \rule{0.5em}{0.5em}}
\begin{document}

\title{Graphs of unitary matrices}
\author{Simone Severini\thanks{Email address:\ ss54@york.ac.uk}\\Department of Computer Science, University of Bristol, U.K.}
\maketitle

\begin{abstract}
The \emph{support }of a matrix $M$ is the $\left(  0,1\right)  $-matrix with
$ij$-th entry equal to $1$ if the $ij$-th entry of $M$ is non-zero, and equal
to $0$, otherwise. The digraph whose adjacency matrix is the support of $M$ is
said to be the \emph{digraph of }$M$. In this paper, we observe some general
properties of digraphs of unitary matrices.


\end{abstract}

\section{Introduction}

A (finite) \emph{directed graph}, for short \emph{digraph}, consists of a
non-empty finite set of elements called \emph{vertices} and a (possibly
empty)\ finite set of ordered pairs of vertices called \emph{arcs}. Let
$D=\left(  V,A\right)  $ be a digraph with vertex-set $V\left(  D\right)  $
and arc-set $A\left(  D\right)  $. In a digraph a \emph{loop} is an arc of the
form $\left(  v_{i},v_{i}\right)  $. In a digraph $D$, if $\left(  v_{i}%
,v_{j}\right)  ,\left(  v_{j},v_{i}\right)  \in A\left(  D\right)  $ the pair
$\left\{  \left(  v_{i},v_{j}\right)  ,\left(  v_{j},v_{i}\right)  \right\}  $
is called \emph{edge} and denoted simply by $\left\{  v_{i},v_{j}\right\}  $.
A digraph $D$ is \emph{symmetric} if, for every $v_{i},v_{j}\in V\left(
D\right)  $, $\left(  v_{i},v_{j}\right)  \in A\left(  D\right)  $ if and only
if $\left(  v_{j},v_{i}\right)  \in A\left(  D\right)  $. Naturally, a
symmetric digraph is also called \emph{graph}. The \emph{adjacency matrix} of
a digraph $D$ on $n$ vertices, denoted by $M\left(  D\right)  $, is the
$n\times n$ $\left(  0,1\right)  $-matrix with $ij$-th entry
\[
M_{i,j}\left(  D\right)  =\left\{
\begin{tabular}
[c]{ll}%
$1$ & if $\left(  v_{i},v_{j}\right)  \in A\left(  D\right)  ,$\\
$0$ & otherwise.
\end{tabular}
\ \right.
\]
Let $M$ be an $n\times n$ matrix (over any field). The \emph{support of }$M$
is the $n\times n$ $\left(  0,1\right)  $-matrix with $ij$-th entry equal to 1
if $M_{i,j}\neq0$, and equal to $0$, otherwise. The \emph{digraph of} $M$ is
the digraph whose adjacency matrix is the support of $M$. If a digraph $D$ is
the digraph of a matrix $M$ then we say that $D$ (or, indistinctly, $M\left(
D\right)  $) \emph{supports} $M$. An $n\times n$ complex matrix $U$ is
\emph{unitary} if $U^{\dagger}U=UU^{\dagger}=I_{n}$, where $U^{\dagger}$ is
the adjoint of $U$ and $I_{n}$ the identity matrix of size $n$. We denote by
$\mathcal{U}$ the set of the digraphs of unitary matrices. Properties of
digraphs of unitary matrices are investigated in \cite{BBS93}, \cite{CJLP99},
\cite{CS00}, \cite{GZ98} and \cite{S03} (see also the references contained
therein). These articles mainly study the number of non-zero entries in
unitary matrices with specific combinatorial properties (\emph{e.g.}
irreducibility, first column(row) without zero-entries, \emph{etc.}). The
present paper observes some structural properties of digraphs and Cayley
digraphs, of unitary matrices. The next two subsections outline the paper.

\subsection{Cayley digraphs}

Section 2 is dedicated to Cayley digraphs. Let $G$ be a finite group and let
$S\subset G$. We denote by $e$ the identity element of a group $G$. We write
$G=\left\langle S:\mathcal{R}\right\rangle $ to mean that $G$ is generated by
$S$ with a set of relations $\mathcal{R}$. When we do not need to specify
$\mathcal{R}$, we write simply $G=\left\langle S\right\rangle $.

The \emph{Cayley digraph of }$G$\emph{\ with respect to }$S$, denoted by
$X\left(  G;S\right)  $, is the digraph whose vertex-set is $G$, and whose
arc-set is the set of all ordered pairs $\{\left(  g,sg\right)  :g\in G,s\in
S\}$. Let $\rho_{reg}$ be the regular permutation representation of $G$. Then
\[%
\begin{tabular}
[c]{lll}%
$M\left(  X\left(  G;S\right)  \right)  =\sum_{i=1}^{k}\rho_{reg}\left(
s_{i}\right)  ,$ & with & $S=\left\{  s_{1},s_{2},...,s_{k}\right\}  .$%
\end{tabular}
\]
Note that $M\left(  X\left(  G;S\right)  \right)  =\widehat{\delta_{S}}\left(
\rho_{reg}\right)  $, where $\widehat{\delta_{S}}\left(  \rho_{reg}\right)  $
is the Fourier transform at $\rho_{reg}$ of the characteristic function of
$S$. For every $S\subset V\left(  D\right)  $, let $N^{-}\left(  S\right)
=\left\{  v_{i}:\left(  v_{i},v_{j}\right)  \in A\left(  D\right)  ,v_{j}\in
S\right\}  $ and $N^{+}\left(  S\right)  =\left\{  v_{j}:\left(  v_{i}%
,v_{j}\right)  \in A\left(  D\right)  ,v_{i}\in S\right\}  $ be the
\emph{in-neighbourhood} and the \emph{out-neighbourhood} of $S$, respectively.
If $D$ is a graph the \emph{neighbourhood} of $S$ is denoted simply by
$N\left(  S\right)  $. A digraph $D$ is $d$\emph{-regular} if, for every
$v_{i}\in V$, $\left\vert N^{-}\left(  v_{i}\right)  \right\vert =\left\vert
N^{+}\left(  v_{i}\right)  \right\vert =d$. A Cayley digraph $X\left(
G;S\right)  $ is on $n=\left\vert G\right\vert $ vertices and $d$-regular,
where $d=\left\vert S\right\vert $. If $S=S^{-1}$ then the Cayley digraph
$X\left(  G;S\right)  $ is called \emph{Cayley graph}. In Section 2, we prove
the following theorem, and construct some examples.

\begin{theorem}
\label{primo}Let $G$ be a group with a generating set of two elements. Then
there exists a set $S\subset G$, such that $G=\left\langle S\right\rangle $
and the Cayley digraph $X\left(  G;S\right)  \in\mathcal{U}$.
\end{theorem}

Since every finite simple group has a generating set of two elements
\cite{AG84}, we have this corollary:

\begin{corollary}
Let $G$ be a finite simple group. Then there exists a set $S\subset G$, such
that $G=\left\langle S\right\rangle $ and the Cayley digraph $X\left(
G;S\right)  \in\mathcal{U}$.
\end{corollary}

Let $\Pi_{n}$ be the group of all permutation matrices of size $n$. Let
$P,Q\in\Pi_{n}$. We say that $P$ and $Q$ are \emph{complementary} if, for any
$1\leq h,i,j,k\leq n$,
\[%
\begin{tabular}
[c]{lll}%
$P_{i,j}=P_{h,k}=Q_{i,k}=1\text{ }$ & imply & $Q_{h,j}=1,$%
\end{tabular}
\]
and, consequently,
\[%
\begin{tabular}
[c]{lll}%
$Q_{i,j}=Q_{h,k}=P_{i,k}=1$ & imply & $P_{h,j}=1.$%
\end{tabular}
\]
We make some observations about Cayley digraphs whose adjacency matrix is sum
of complementary permutations. We show that the $n$-cube is the digraph of a
unitary matrix. This is also true for the de Bruijn digraph \cite{ST}. It
might be interesting to remark that the $n$-cube and the de Bruijn digraphs
are among the best-known architectures for interconnection networks (see,
\emph{e.g.}, \cite{H97}, for a survey of this subject, with particular
attention to Cayley digraphs). It would be interesting to deepen the study of
digraphs of unitary matrices seen as specific architectures for
interconnection networks.

\subsection{Digraphs}

Section 3 is dedicated to digraphs in general. Theorem 3 is the main result of
the section. A \emph{dipath} is a non-empty digraph $D$, where $V\left(
D\right)  =\left\{  v_{0},v_{1},...,v_{k}\right\}  $ and $A\left(  D\right)
=\{\left(  v_{0},v_{1}\right)  ,\left(  v_{1},v_{2}\right)  ,...,\left(
v_{k-1},v_{k}\right)  \}$. The vertex $v_{0}$ is the \emph{initial vertex} of
$D$; the vertex $v_{k}$ is the \emph{final vertex}. We say that $D$ is a
\emph{dipath from }$v_{0}$\emph{\ to }$v_{k}$. A \emph{path} is a non-empty
graph $D$, with $V\left(  D\right)  =\left\{  v_{0},v_{1},...,v_{k}\right\}  $
and $A\left(  D\right)  =\left\{  \left\{  v_{0},v_{1}\right\}  ,\left\{
v_{1},v_{2}\right\}  ,...,\left\{  v_{k-1},v_{k}\right\}  \right\}  $. Two or
more dipaths (paths) are \emph{independent} if none of them contains an inner
vertex of another. A digraph is \emph{connected} if, for every $v_{i}$ and
$v_{j}$, there is a dipath from $v_{i}$ to $v_{j}$, or \emph{viceversa};
\emph{strongly-connected} if, for every $v_{i}$ and $v_{j}$, there is a dipath
from $v_{i}$ to $v_{j}$ and to $v_{j}$ to $v_{i}$. A Cayley digraph is
strongly-connected. A $k$\emph{-dicycle} is a dipath on $k$ arcs in which the
initial and final vertex coincide. If all the vertices and the arcs of a
dipath (dicycle) are \emph{all} distinct then the dipath (dicycle) is an
\emph{Hamilton dipath} (\emph{dicycle}). A digraph\ spanned by an Hamilton
dicycle is said to be \emph{hamiltonian}. In a graph, the analogue of dicycle
and hamiltonian dipath (dicycle) are called \emph{cycle} and \emph{hamiltonian
path} (\emph{cycle}). In a digraph $D$ on $n\geq2$ vertices, a
\emph{disconnecting set} of arcs (edges)\ is a subset $T\subset A\left(
D\right)  $ such that $D-T$ has more connected components than $D$. The
\emph{arc(edge)-connectivity} is the smallest number of edges in any
disconnecting set. A \emph{cut} of $D$ is a subset $S\subset V\left(
D\right)  $ such that $D-S$ has more connected components than $D$. The
\emph{vertex-connectivity} of $D$ is the smallest number of vertices in any
cut of $D$. A digraph $D$ is said to be $k$\emph{-vertex-connected}
($k$\emph{-arc(edge)-connected}) if its vertex-connectivity
(arc(edge)-connectivity) is larger or equal than $k$. A \emph{cut-vertex}, a
\emph{directed bridge}, and a \emph{bridge}, are respectively a vertex, an
arc, and an edge, whose deletion increases the number of connected components
of $D$. A digraph is \emph{inseparable} if it is without cut-vertices;
\emph{bridgeless} if it is without bridges. Let $\overleftrightarrow{K}_{2}$
and $\overleftrightarrow{K}_{2}^{+}$be respectively the complete graph on two
vertices and the complete graph on two vertices with a self-loop at each
vertex. We prove the following theorem, and state some of its natural corollaries.

\begin{theorem}
\label{th1}Let $D$ be a digraph. If $D\in\mathcal{U}$ then:

\begin{enumerate}
\item $D$ is without directed bridges;

\item $D$ is bridgeless, unless the bridge is in a connected component that is
either $\overleftrightarrow{K}_{2}$ or $\overleftrightarrow{K}_{2}^{+}$.

\item $D$ is inseparable, unless a cut-vertex is in a connected component that
is either $\overleftrightarrow{K}_{2}$ or $\overleftrightarrow{K}_{2}^{+}$.
\end{enumerate}
\end{theorem}

\subsection{A motivation}

Unitary matrices appear in many areas of Physics and are of fundamental
importance in Quantum Mechanics. The \emph{time-evolution} of the state of an
$n$-level quantum system, assumed to be isolated from the environment, is
reversible and determined by the rubric $\rho\longrightarrow U_{t}\rho
U_{t}^{-1}$, where $\left\{  U_{t}:-\infty<t<\infty\right\}  $ is a continuous
group of unitary matrices, and $\rho$, the \emph{state of the system}, is an
$n\times n$ Hermitian matrix, which is positive definite and has unit trace.
Sometimes it is useful to look at a quantum system as evolving
\emph{discretely}, under the same unitary matrix: $\rho\longrightarrow U\rho
U^{-1}\longrightarrow U^{2}\rho U^{-2}\longrightarrow\cdots\longrightarrow
U^{n}\rho U^{-n}$. Suppose that to an $n$-level quantum system is assigned a
digraph $D$ on $n$ vertices, in the following sense: the vertices of $D$ are
labeled by given states of the system; the arc $\left(  v_{i},v_{j}\right)  $
means non-zero probability of transition from the state labeled by $v_{i}$ to
the state labeled by $v_{j}$, in one time-step, that is in one application of
$U$. As it happens for random walks on graphs, important features of this
evolution depend on the combinatorial properties of $D$, the digraph of the
\textquotedblleft transition matrix\textquotedblright\ $U$ (here $U$ unitary
rather than stochastic). Quantum evolution in digraphs have recently drawn
attention in Quantum Computation (see \emph{e.g.}, \cite{AAKV01}, \cite{SKW02}
and \cite{C+03}) and in the study of statistical properties of quantum systems
in relation to Random Matrix Theory (see, \emph{e.g.}, \cite{KS99},
\cite{T01}, \cite{KS03}, \cite{ST} and the references contained therein).

\section{Cayley digraphs}

\subsection{Proof of Theorem 1}

The \emph{line digraph} of a digraph $D$, denoted by $\overrightarrow{L}D$, is
defined as follows:\ the vertex-set of $\overrightarrow{L}D$ is $A\left(
D\right)  $; $\left(  v_{i},v_{j}\right)  ,\left(  v_{k},v_{l}\right)  \in
A\left(  \overrightarrow{L}D\right)  $ if and only if $v_{j}=v_{k}$. The
digraph $D$ is said to be the \emph{base} of $\overrightarrow{L}D$. (See,
\emph{e.g.}, \cite{P96}, for a survey on line graphs and digraphs.)

\begin{definition}
[Independent full submatrix]A rectangular array, say $M^{\prime}$, of entries
from an $n\times n$ matrix $M$ is an\emph{\ independent full submatrix} when,
if $M_{i,j}\in M^{\prime}$ then, for every $1\leq k,l\leq n$, either
$M_{i,k}\in M^{\prime}$ or $M_{i,k}=0$, and, either $M_{l,j}\in M^{\prime}$ or
$M_{l,j}=0$. In addition, if $M_{i,j}\in M^{\prime}$ then $M_{i,j}\notin
M^{\prime\prime}$, where $M^{\prime\prime}$ is an independent full submatrix
different from $M^{\prime}$.
\end{definition}

\begin{example}
Consider the matrix
\[
M=\left[
\begin{array}
[c]{ccccc}%
0 & 0 & x_{1,3} & x_{1,4} & x_{1,5}\\
0 & 0 & x_{2,3} & x_{2,4} & x_{2,5}\\
x_{3,1} & x_{3,2} & 0 & 0 & 0\\
x_{4,1} & x_{4,2} & 0 & 0 & 0\\
0 & 0 & x_{5,3} & x_{5,4} & x_{5,5}%
\end{array}
\right]  .
\]
The matrices
\[%
\begin{tabular}
[c]{lll}%
$\left[
\begin{array}
[c]{cc}%
x_{3,1} & x_{3,2}\\
x_{4,1} & x_{4,2}%
\end{array}
\right]  $ & and & $\left[
\begin{array}
[c]{ccc}%
x_{1,3} & x_{1,4} & x_{1,5}\\
x_{2,3} & x_{2,4} & x_{2,5}\\
x_{5,3} & x_{5,4} & x_{5,5}%
\end{array}
\right]  $%
\end{tabular}
\]
are independent full submatrices of $M$.
\end{example}

The following is an easy lemma. This can be also seen as a corollary of
Theorem 2.15 in \cite{S03}.

\begin{lemma}
\label{regular}Let $D$ be a Cayley digraph. If there exists a digraph
$D^{\prime}$ such that $D=\overrightarrow{L}D^{\prime}$ then $D\in\mathcal{U}$.
\end{lemma}

\begin{proof}
By the \emph{Richard characterization of line digraphs} (see, \emph{e.g.},
\cite{P96}), a digraph $D$ is a line digraph if and only if the following two
conditions hold:

\begin{itemize}
\item The columns of $M\left(  D\right)  $ are identical or orthogonal.

\item The rows of $M\left(  D\right)  $ are identical or orthogonal.
\end{itemize}

This means that, if $D$ is a line digraph then every non-zero entry of
$M\left(  D\right)  $ belongs to an independent full submatrix. Moreover if
$D$ is a regular line digraph then \emph{all} the independent full submatrices
of $M\left(  D\right)  $ are square. Suppose that $D$ is a Cayley digraph and
a line digraph. Observe that:

\begin{itemize}
\item[(i)] Since $D$ is strongly-connected, $M\left(  D\right)  $ has neither
zero-rows nor zero-columns.

\item[(ii)] Since $D$ is regular, every independent full submatrix of
$M\left(  D\right)  $ is square.
\end{itemize}

Combining (i) and (ii), and since the all-ones matrix supports a unitary
matrix, the lemma follows.
\end{proof}

Let $\mathbb{Z}_{n}$\ be the additive group of the integers modulo $n$.

\begin{remark}
\label{rem}The converse of Lemma \ref{regular} is false. For example, consider
the Cayley digraph $D=X\left(  \mathbb{Z}_{4};\left\{  1,2,3\right\}  \right)
$. The adjacency matrix of $D$ is
\[
M\left(  D\right)  =\left[
\begin{array}
[c]{cccc}%
0 & 1 & 1 & 1\\
1 & 0 & 1 & 1\\
1 & 1 & 0 & 1\\
1 & 1 & 1 & 0
\end{array}
\right]  \emph{.}%
\]
The matrix
\[
U=\frac{1}{^{\sqrt{3}}}\left[
\begin{array}
[c]{cccc}%
0 & 1 & 1 & 1\\
1 & 0 & -1 & 1\\
1 & 1 & 0 & -1\\
1 & -1 & 1 & 0
\end{array}
\right]
\]
is unitary. Since $M\left(  D\right)  $ supports $U$, $D\in\mathcal{U}$. Note
that $D$ is not a line digraph since it does not satisfy the Richard characterization.
\end{remark}

A \emph{multidigraph} is a digraph with possibly more than one arc $\left(
v_{i},v_{j}\right)  $, for some $v_{i}$ and $v_{j}$.

\begin{lemma}
[Mansilla-Serra, \cite{MS01}]\label{sonia}Let $G=\left\langle S\right\rangle $
be a finite group. If, for some $x\in S^{-1}$, $xS=H$, where $H$ is a subgroup
of $G$ such that $\left|  H\right|  =k=\left|  S\right|  $, then $X\left(
G;S\right)  =\overrightarrow{L}D$, where $D$ is a $k$-regular (multi)digraph.
\end{lemma}

Let $C_{n}$ be a cyclic group of order $n$. The proof of Theorem \ref{primo}
makes use of Lemma \ref{regular} and Lemma \ref{sonia}.

\bigskip

\begin{proof}
[\textbf{Proof of Theorem 1}]Let $G=\left\langle S\right\rangle $, where
$S=\left\{  s_{1},s_{2}\right\}  $. Take $s_{1}^{-1}\in S^{-1}$ (or,
equivalently, $s_{2}^{-1}$). Then $s_{1}^{-1}S=\left\{  s_{1}^{-1}s_{1}%
,s_{1}^{-1}s_{2}\right\}  =\left\{  e,s_{1}^{-1}s_{2}\right\}  $. Let
$s_{1}^{-1}s_{2}$ have order $n$. Consider $C_{n}=\left\langle s_{1}^{-1}%
s_{2}\right\rangle $. Write $T=s_{1}C_{n}$. Then $C_{n}=s_{1}^{-1}T$. Now,
observe that $s_{1}s_{1}^{-1}s_{2}=s_{2}$ and $s_{1}\left(  s_{1}^{-1}%
s_{2}\right)  ^{-1}s_{1}^{-1}s_{2}=s_{1}s_{2}^{-1}s_{1}s_{1}^{-1}s_{2}=s_{1}$.
Then $S\subset T$ and $G=\left\langle T\right\rangle $. Since $T$ is a left
coset of $C_{n}$, $\left\vert T\right\vert =n=\left\vert C_{n}\right\vert $.
By Lemma \ref{sonia}, $X\left(  G;T\right)  $ is a line digraph, and hence, by
Lemma \ref{regular}, $X\left(  G;T\right)  \in\mathcal{U}$.
\end{proof}

\subsection{Examples}

Let $D_{n}$\ be a dihedral group of order $2n$.

\begin{example}
The standard presentation of $D_{n}$\ (see, e.g., \cite{CM72}, \S 1.5) is
\[
\left\langle s_{1},s_{2}:s_{1}^{n}=s_{2}^{2}=e,s_{2}s_{1}s_{2}=s_{1}%
^{-1}\right\rangle .
\]
By Lemma \ref{regular} $X\left(  D_{n};\left\{  s_{1},s_{2}\right\}  \right)
\in\mathcal{U}$\ since (see, \emph{e.g.}, \cite{BEFS95}),
\[
X\left(  D_{n},\left\{  s_{1},s_{2}\right\}  \right)  \cong\overrightarrow
{L}X\left(  \mathbb{Z}_{n},\left\{  1,n-1\right\}  \right)
\]

\end{example}

\begin{definition}
[Digraph $P\left(  n,k\right)  $, \cite{F84}]Given integers $k$ and $n$,
$1\leq k\leq n-1$, $P\left(  n,k\right)  $ is the digraph whose vertices are
the permutations on $k$-tuples from the set $\left\{  1,2,...,n\right\}  $ and
whose arcs are of the form $\left(  \left(  i_{1}\text{ }i_{2}\text{ ...
}i_{k}\right)  ,\left(  i_{2}\text{ }i_{3}\text{ ... }i_{k}i\right)  \right)
$, where $i\neq i_{1},i_{2},...,i_{k}$.
\end{definition}

Let $S_{n}$\ be a symmetric group on a set $\left\{  1,2,...,n\right\}  $.

\begin{example}
Let $S_{n}=\left\langle s_{1},s_{2}\right\rangle $, where $s_{1}=\left(
1\text{ }2\text{ ... }n\right)  $ and $s_{2}=\left(  1\text{ }2\text{ ...
}n-1\right)  $\emph{. }This is because (see, \emph{e.g.}, \cite{CM72}, \S 1.7)
$S_{n}=\left\langle \left(  1\text{ }2\text{ ... }n\right)  ,\left(  1\text{
}n\right)  \right\rangle $ and $\left(  1\text{ }n\right)  =\left(  1\text{
}2\text{ ... }n-1\right)  \cdot\left(  1\text{ }2\text{ ... }n\right)  ^{-1}$.
By Lemma \ref{regular} and since (\cite{BFF97}, Lemma 2.1) $X\left(
S_{n};\left\{  s_{1},s_{2}\right\}  \right)  \cong\overrightarrow{L}P\left(
n,n-2\right)  $, we have $X\left(  S_{n};\left\{  s_{1},s_{2}\right\}
\right)  \in\mathcal{U}$.
\end{example}

\begin{example}
The Cayley digraph $X\left(  S_{n};T\right)  $, where $T=\left(  1\text{
}2\right)  C_{n-1}$, is the digraph of a unitary matrix. Consider
$S_{n}=\left\langle S=\left\{  \left(  1\text{ }2\right)  ,\left(  1\text{
}2\text{ ... }n\right)  \right\}  \right\rangle $, Then $S^{-1}=\left\{
\left(  1\text{ }2\right)  ,\left(  1\text{ }n\text{ ... }2\right)  \right\}
$. Write $x=\left(  1\text{ }2\right)  \in S^{-1}$. Then
\[
\left(  1\text{ }2\right)  S=\left\{  e,\left(  2\text{ }3\text{ ...
}n\right)  \right\}  \emph{.\ }%
\]
Consider $C_{n-1}=\left\langle e,\left(  2\text{ }3\text{ ... }n\right)
\right\rangle $. Write
\[
T=\left(  1\text{ }2\right)  C_{n-1}=\left\{  \left(  1\text{ }2\right)
,\left(  1\text{ }2\right)  \left(  2\text{ }3\text{ ... }n\right)  =\left(
1\text{ }2\text{ ... }n\right)  ,...\right\}  .
\]
Since $S\subset T$, $S_{n}=\left\langle T\right\rangle $. Moreover $x\in
T^{-1}$, $C_{n-1}=xT\ $and $\left\vert T\right\vert =n-1=C_{n-1}$. Then, by
Lemma \ref{regular} and Lemma \ref{sonia}, $X\left(  S_{n};T\right)
\in\mathcal{U}$.
\end{example}

\subsection{Cayley digraphs of abelian groups}

\subsubsection{General properties}

Let $conv\left\{  P_{1},...,P_{m}\right\}  $ be the convex hull of the
matrices $P_{1},...,P_{m}\in\Pi_{n}$. Note that all the matrices that belong
this convex hull have the same digraph. A \emph{doubly-stochastic} matrix is a
non-negative matrix whose row sums and column sums give one. The Birkhoff
theorem for doubly-stochastic matrices (see, \emph{e.g.}, \cite{B97}) says
that the set of $n\times n$ doubly-stochastic matrices\ is the convex hull of
permutation matrices. A doubly-stochastic matrix $M$\ is \emph{uni-stochastic}%
\ if $M_{i,j}=\left|  U_{i,j}\right|  ^{2}$. The existence of a
``Birkhoff-type'' theorem for uni-stochastic matrices is an open problem (see,
\emph{e.g.}, \cite{F88} and \cite{L97}).

Let $P_{1},...,P_{m}$, such that $conv\left\{  P_{1},...,P_{m}\right\}
\subset\mathcal{O}$, where $\mathcal{O}$ denotes the set of uni-stochastic
matrices. If a digraph $D$ supports a uni-stochastic matrix then $D$ supports
a unitary matrix, and \emph{viceversa}. If $D$ supports $conv\left\{
P_{1},...,P_{m}\right\}  \subset\mathcal{O}$ then, obviously, $D$ supports a
unitary matrix. In such a case, with an abuse of notation, we write
$D\in\mathcal{O}$.

\begin{theorem}
[Au-Young and Cheng, \cite{AC91}]\label{0}Let $P_{1},...,P_{m}\in\Pi_{n}$. If
\[
conv\left\{  P_{1},...,P_{m}\right\}  \subset\mathcal{O}%
\]
then $P_{1},...,P_{m}$ are pairwise complementary.
\end{theorem}

\begin{proposition}
\label{1}If $X\left(  G;S\right)  \in\mathcal{O}$ then:

\begin{enumerate}
\item For every $s,t\in S$, and $1\leq h,i,j,k\leq\left|  G\right|  $, if
$g_{j}=sg_{i},g_{k}=sg_{h}$ and $g_{k}=tg_{i}$ then $g_{j}=tg_{h}$.

\item For every $s,t\in S$, $st^{-1}=ts^{-1}$.

\item The order of $G$ is even.

\item If $G$ is abelian then, for every $s,t\in S$, $2s=2t$.
\end{enumerate}
\end{proposition}

\begin{proof}
In the order:

\begin{enumerate}
\item By Theorem \ref{0}.

\item From 1, since $s^{-1}g_{j}=g_{i}$ and $t^{-1}g_{k}=g_{i}$, we have
$s^{-1}g_{j}=t^{-1}g_{k}$, and, since $s^{-1}g_{k}=g_{h}$ and $t^{-1}%
g_{j}=g_{h}$, we have $s^{-1}g_{k}=t^{-1}g_{h}$. Then, since $g_{k}%
=st^{-1}g_{j}$, we obtain $s^{-1}g_{j}=t^{-1}st^{-1}g_{j}$ and $ts^{-1}%
g_{j}=st^{-1}g_{j}$, that implies $st^{-1}=ts^{-1}$. Since $st^{-1}%
=ts^{-1}=\left(  st^{-1}\right)  ^{-1}$, $st^{-1}$ is an involution.

\item From point 2, since a group of odd order is without involutions.

\item From point 2, since $G$ abelian, given that $s=ts^{-1}t=s^{-1}2t$, we
have $2s=2t$.
\end{enumerate}
\end{proof}

\subsubsection{Cayley digraphs of cyclic groups}

\begin{proposition}
\label{2}If $X\left(  \mathbb{Z}_{n};S\right)  \in\mathcal{O}$ then:

\begin{enumerate}
\item $\left\vert S\right\vert =2$ and $t=s+\frac{n}{2}\left(
\operatorname{mod}n\right)  $.

\item $\mathbb{Z}_{n}=\left\langle s,t\right\rangle $ if and only if $s$ odd,
or $s$ even and $n=4m+2$ (that is $\frac{n}{2}$ is odd) where $m$ is a
non-negative integer.

\item If $\mathbb{Z}_{n}=\left\langle s,t\right\rangle $ then $X\left(
\mathbb{Z}_{n};\left\{  s,t\right\}  \right)  $ is hamiltonian.

\item $X\left(  \mathbb{Z}_{n};\left\{  s,t\right\}  \right)  $ is a graph if
and only if $s=\frac{n}{4}$.

\item If $X\left(  \mathbb{Z}_{n};\left\{  s,t\right\}  \right)  $ is a graph
and $\mathbb{Z}_{n}=\left\langle s,t\right\rangle $ then it is not hamiltonian.

\item $\left\vert N^{+}\left(  s\right)  \cap N^{+}\left(  t\right)
\right\vert =\left\vert N^{-}\left(  s\right)  \cap N^{-}\left(  t\right)
\right\vert =2$.
\end{enumerate}
\end{proposition}

\begin{proof}
In the order:

\begin{enumerate}
\item From 4 of Proposition \ref{1}, $2s=2t$. Let $t>s$. Then $t=s+x\left(
\operatorname{mod}n\right)  $ and $2s+2x\left(  \operatorname{mod}n\right)
=2t\left(  \operatorname{mod}n\right)  $. This occurs if and only if
$x=\frac{n}{2}$. Then $t=s+\frac{n}{2}\left(  \operatorname{mod}n\right)  $.

\item If $s$ and $t$ are both even then they generate the even subgroup of
order $n/2$. In the other cases, $s$ and $t$ generate $\mathbb{Z}_{n}$.
Clearly if $s$ even then $t$ odd if and only if $n/2$ is odd, that is $n=4m+2$.

\item From point 2, either $s$ or $t$ has to be odd. Since $n$ is even, an odd
element of $\mathbb{Z}_{n}$ has order $n$. Suppose $s$ odd. In $X\left(
\mathbb{Z}_{n};\left\{  s,t\right\}  \right)  $ there is then an hamiltonian
cycle $e,s,2s,...,\left(  n-1\right)  s,e$.

\item From point 1, $t=s+\frac{n}{2}\left(  \operatorname{mod}n\right)  $. If
$X\left(  \mathbb{Z}_{n};\left\{  s,t\right\}  \right)  $ is a graph then
$t=s^{-1}$, that is $t=n-s\left(  \operatorname{mod}n\right)  $. So,
$n-s\left(  \operatorname{mod}n\right)  =s+\frac{n}{2}\left(
\operatorname{mod}n\right)  $, which implies $s=\frac{n}{4}$. The sufficiency
is clear.

\item From the previous point, $X\left(  \mathbb{Z}_{n};\left\{  s,t\right\}
\right)  $ is a graph if and only if $s=\frac{n}{4}$. Then $n=4m$. From 3,
since $s$ is even, we need $t$ odd to generated $\mathbb{Z}_{n}$. From 2, $t$
is odd if $n=4m+2$. A contradiction.
\end{enumerate}
\end{proof}

The \emph{distance} from a vertex $v_{i}$ to a vertex $v_{j}$ is denoted by
$d\left(  v_{i},v_{j}\right)  $ and it is the length (the number of arcs) of
the shortest dipath from $v_{i}$ to $v_{j}$. The \emph{diameter} of $D=\left(
V,A\right)  $ is $dia\left(  D\right)  =\max_{\left(  v_{i},v_{j}\right)  \in
V\times V}d\left(  v_{i},v_{j}\right)  $.

\begin{proposition}
\label{mult}If $X\left(  \mathbb{Z}_{n};\left\{  s,t\right\}  \right)
\in\mathcal{O}$ then:

\begin{enumerate}
\item It is a line digraph of the multidigraph with adjacency matrix $M=2\cdot
M\left(  X\left(  \mathbb{Z}_{n/2};\left\{  1\right\}  \right)  \right)  $.

\item $dia\left(  X\left(  \mathbb{Z}_{n};\left\{  s,t\right\}  \right)
\right)  =\frac{n}{2}+1$.
\end{enumerate}
\end{proposition}

\begin{proof}
In the order:

\begin{enumerate}
\item From 6 of Proposition \ref{2}, follows that the rows and columns of
$M\left(  D\right)  $ are identical or orthogonal. By the Richard
characterization (cfr. proof of Lemma \ref{regular}), this is sufficient for a
digraph to be a line digraph. Observe that the base digraph of $X\left(
\mathbb{Z}_{n};\left\{  s,t\right\}  \right)  $ is the multidigraph with
adjacency matrix $M.$

\item Since $dia\left(  X\left(  \mathbb{Z}_{n/2};\left\{  1\right\}  \right)
\right)  =n/2$ and since the diameter of the line digraph increases of one
unit in respect to the diameter of its base digraph (see, \emph{e.g.},
\cite{P96}), the proposition follows.
\end{enumerate}
\end{proof}

\begin{remark}
Let $D=X\left(  \mathbb{Z}_{n};\left\{  s,t\right\}  \right)  \in\mathcal{O}$.
From Proposition \ref{mult}, follows that the eigenvalues of $D$ are the $n/2$
eigenvalues of the multidigraph, which are $\left\{  2\omega^{j}:0\leq
j<n/2,\omega=\allowbreak e^{4i\frac{\pi}{n}}\right\}  $, plus an eigenvalue
zero with multiplicity $n/2$.
\end{remark}

An \emph{automorphism} of a digraph $D$ is a permutation $\pi$ of $V\left(
D\right)  $, such that $\left(  v_{i},v_{j}\right)  \in A\left(  D\right)  $
if and only if $\left(  \pi\left(  v_{i}\right)  ,\pi\left(  v_{j}\right)
\right)  \in A\left(  D\right)  $. Let $Aut\left(  D\right)  $ be the group of
the automorphisms of a digraph $D$. It is well-known that if $D=X\left(
G;S\right)  $ is a Cayley digraph then $Aut\left(  D\right)  $ contains the
regular representation of $G$. This implies that a Cayley digraph is
\emph{vertex-transitive}, that is its automorphism group acts transitively on
its vertex-set. A digraph $D$ is \emph{arc-transitive} if, for any pair of
arcs $\left(  v_{i},v_{j}\right)  $ and $\left(  v_{k},v_{l}\right)  $, there
exists a permutation $\pi\in Aut\left(  D\right)  $ such that $\pi\left(
v_{i}\right)  =v_{k}$ and $\pi\left(  v_{j}\right)  =v_{l}$.

\begin{proposition}
If $X\left(  \mathbb{Z}_{n};\left\{  s,t\right\}  \right)  \in\mathcal{O}$
then it is arc-transitive.
\end{proposition}

\begin{proof}
Let $D=X\left(  \mathbb{Z}_{n};\left\{  s,t\right\}  \right)  $. Since $D$ is
a Cayley digraph, $Aut\left(  D\right)  $ contains the regular representation
of $D$. Take the element $\frac{n}{2}$ and look at $\frac{n}{2}$ has an
automorphism of $D$. The action of $\frac{n}{2}$ on $s$ gives $s+\frac{n}%
{2}\left(  \operatorname{mod}n\right)  $. The action of $\frac{n}{2}$ on
$s+\frac{n}{2}\left(  \operatorname{mod}n\right)  $ gives $s$. The proposition
follows easily, $\frac{n}{2}\left(  S\right)  =S$, and $\frac{n}{2}$ can be
seen as a group homomorphism.
\end{proof}

\begin{remark}
Consider a nearest neighbor random walk on $X\left(  \mathbb{Z}_{n};\left\{
s,t\right\}  \right)  \in\mathcal{O}$, with probability $p\left(  s\right)
=\frac{1}{2}=q\left(  t\right)  $. This random walk is non-ergodic since
$\gcd\left(  t-s,n\right)  =\frac{n}{2}$. Observe that, in Cesaro-mean, the
random walk is ergodic and converges in $n/2$ steps towards uniformity. It
would be interesting to observe if random walks on digraphs of unitary
matrices have a characteristic behaviour.
\end{remark}

\subsubsection{General abelian groups}

Let $G=\mathbb{Z}_{p_{1}}\times\mathbb{Z}_{p_{2}}\times\cdots\times\mathbb{Z}%
$, be an abelian group written in its prime-power canonical form. An element
of $G$ has then the form $\left(  g_{1},g_{2},...,g_{l}\right)  $. Let
\[
S=\left\{  \left(  s_{1_{1}},s_{2_{1}},...,s_{l_{1}}\right)  ,\left(
s_{1_{2}},s_{2_{2}},...,s_{l_{2}}\right)  ,...,\left(  s_{1_{k}},s_{2_{k}%
},...,s_{l_{k}}\right)  \right\}
\]
be a set of generators of $G$. If $X\left(  G;S\right)  \in\mathcal{O}$ then,
from 4 of Proposition \ref{1}, $2s=2t$, for every $s,t\in S$. Then, for every
$i$ and $j$,
\begin{align*}
2\left(  s_{1_{i}},s_{2_{i}},...,s_{l_{i}}\right)   &  =2\left(  s_{1_{j}%
},s_{2_{j}},...,s_{l_{j}}\right)  =\\
\left(  2s_{1_{i}},2s_{2_{i}},...,2s_{l_{i}}\right)   &  =\left(  2s_{1_{j}%
},2s_{2_{j}},...,2s_{l_{j}}\right)  .
\end{align*}

\begin{proposition}
Let $G$ be abelian and let $X\left(  G;S\right)  \in\mathcal{O}$.

\begin{enumerate}
\item If $p_{i}$ is odd then $s_{i_{j}}=s_{i_{k}}$, for every $j$ and $k$.

\item If every $p_{i}$ is odd then $\left|  S\right|  =1$.
\end{enumerate}
\end{proposition}

\begin{proof}
In the order:

\begin{enumerate}
\item Suppose that $p_{i}$ is odd. From 4 of Proposition \ref{1}, $2s_{i_{j}%
}=2s_{i_{k}}$. The result follows. This implies that, if $G=\left\langle
S\right\rangle $ then $s_{i_{j}}\neq e$. In fact, if $s_{i_{j}}=e$ then, for
every $k$, $s_{i_{k}}=e$, and, in such a case, $G\neq\left\langle
S\right\rangle $.

\item It is a consequence of the previous point.
\end{enumerate}
\end{proof}

\subsubsection{An example:\ the $n$-cube}

An $n$\emph{-cube} (or, equivalently, $n$\emph{-dimensional hypercube}),
denoted by $Q_{n}$, is a graph whose vertices are the vectors of the
$n$-dimensional vector space over the field $GF\left(  2\right)  $. There is
an edge between two vertices of the $n$-cube whenever their Hamming distance
is exactly 1, where the \emph{Hamming distance} between two vectors is the
number of coordinates in which they differ. The $n$-cube is widely used as
architecture for interconnection networks (see, \emph{e.g.}, \cite{H97}). The
$n$-cube is the Cayley digraph of the group $\mathbb{Z}_{2}^{n}$, generated by
the set $S=\left\{  \left(  1,0,...,0\right)  ,\left(  0,1,0,...,0\right)
,...,\left(  0,...,0,1\right)  \right\}  $. Since, for every $s,t\in S$,
$2s=2t$, we have $X\left(  \mathbb{Z}_{2}^{n};S\right)  \in\mathcal{O}$. We
observe this explicitly. Label the vertices of $Q_{n}$ with the binary numbers
representing $0,1,...,2^{n}-1$. Consider
\[
M\left(  Q_{2}\right)  )=\left[
\begin{array}
[c]{cccc}%
0 & 1 & 1 & 0\\
1 & 0 & 0 & 1\\
1 & 0 & 0 & 1\\
0 & 1 & 1 & 0
\end{array}
\right]  .
\]
A \emph{weighing matrix of size }$n$\emph{\ and weight} $k$, denoted by
$W\left(  k,n\right)  $, is a $\left(  -1,0,1\right)  $-matrix such that
$W\left(  k,n\right)  \cdot W\left(  k,n\right)  ^{\intercal}=kI_{n}$.
Clearly, $\frac{1}{\sqrt{k}}W\left(  k,n\right)  $ is unitary. The matrix
\[
M=\left[
\begin{array}
[c]{cccc}%
0 & -1 & 1 & 0\\
-1 & 0 & 0 & 1\\
1 & 0 & 0 & 1\\
0 & 1 & 1 & 0
\end{array}
\right]
\]
is a symmetric weighing matrix, $W\left(  2,4\right)  $. In fact,
$M=M^{\intercal}$ and $MM^{\intercal}=2I_{4}$. Then $M\left(  Q_{2}\right)  $
supports a unitary matrix. The graph $Q_{n}$ is constructed of two copies of
$Q_{n-1}$, where the corresponding vertices of each subgraph are connected.
The base of the construction is the graph with one vertex. The matrix
\[
W\left(  3,8\right)  =\left[
\begin{array}
[c]{cc}%
W\left(  2,4\right)  & -I_{4}\\
I_{4} & W\left(  2,4\right)
\end{array}
\right]
\]
is supported by $Q_{3}$ and is again a weighing matrix, since $W\left(
2,4\right)  $ is symmetric. Note that $W\left(  3,8\right)  $ is not
symmetric. So, in general,
\[
W\left(  k,2^{k}\right)  =\left[
\begin{array}
[c]{cc}%
W\left(  k-1,2^{k-1}\right)  & -I_{2^{k-1}}\\
I_{2^{k-1}} & W\left(  k-1,2^{k-1}\right)  ^{\intercal}%
\end{array}
\right]  ,
\]
is a weighing matrix supported by $Q_{n}$. Note that if $A$ is an $n\times n$
unitary matrix then the block-matrix
\[
\left[
\begin{array}
[c]{cc}%
A & -I_{n}\\
I_{n} & A^{\intercal}%
\end{array}
\right]
\]
is unitary under renormalization, since
\[
\left[
\begin{array}
[c]{cc}%
A & -I_{n}\\
I_{n} & A^{\intercal}%
\end{array}
\right]  \cdot\left[
\begin{array}
[c]{cc}%
A^{\intercal} & I_{n}\\
-I_{n} & A
\end{array}
\right]  =\left[
\begin{array}
[c]{cc}%
2I_{n} & 0\\
0 & 2I_{n}%
\end{array}
\right]
\]

\begin{remark}
The graph obtained by adding a self-loop at each vertex of the $n$-cube also
supports a unitary matrix. The unitary matrix
\[
\frac{1}{\sqrt{3}}\left[
\begin{array}
[c]{cccc}%
1 & 1 & -1 & 0\\
1 & -1 & 0 & 1\\
-1 & 0 & -1 & 1\\
0 & 1 & 1 & 1
\end{array}
\right]
\]
can be seen as the building-block of the iteration.
\end{remark}

\begin{remark}
Let $D$ be a Cayley digraph. If $D\in\mathcal{U}$ then $D$ is not necessarily
a line digraph. The $n$-cube is a counter-example. In fact, $Q_{n}%
\in\mathcal{U}$ and it is not a line digraph.
\end{remark}

\begin{remark}
Let $D$ be the digraph on $n$ vertices, and let $D\in\mathcal{U}$. Let $U$ be
a unitary matrix supported by $D$. If the order of $U\in U\left(  n\right)  $
is $k$ then $U$ generates a cyclic group $\left\{  U,U^{2},...,U^{k}%
=I_{n}\right\}  \subset U\left(  n\right)  $. It would be interesting to study
if the digraphs of the matrices $U,U^{2},...,U^{k-1}$ have some common
properties, apart from, trivially, the same number of vertices.
\end{remark}

\section{Digraphs in general}

\subsection{Proof of Theorem 2}

\begin{proof}
[\textbf{Proof of Theorem 2}]Suppose that $\left(  v_{i},v_{j}\right)  \in
A\left(  D\right)  $ is a directed bridge. We can then label $V\left(
D\right)  $ such that $M\left(  D\right)  $ has the form
\[
\left[
\begin{array}
[c]{ccccccc}
&  &  & 0 &  &  & \\
& M_{1} &  & \vdots &  & \mathbf{0} & \\
&  &  & 0 &  &  & \\
M_{i,1} & \cdots & M_{i,i} & 1 & 0 & \cdots & 0\\
&  &  & M_{i+1,j} &  &  & \\
& \mathbf{0} &  & \vdots &  & M_{2} & \\
&  &  & M_{n,j} &  &  &
\end{array}
\right]  ,
\]
where $M_{i,1},...,M_{i,i},M_{i+1,j},...,M_{n,j}\in\left\{  0,1\right\}  $.
Suppose that $D\in\mathcal{U}$. Then $D$ is \emph{quadrangular} (the term
\textquotedblleft quadrangular\textquotedblright\ has been coined in
\cite{GZ98}), that is, for every $v_{i}$ and $v_{j}$, $\left\vert N^{-}\left(
v_{i}\right)  \cap N^{-}\left(  v_{j}\right)  \right\vert \neq1$ and
$\left\vert N^{+}\left(  v_{i}\right)  \cap N^{+}\left(  v_{j}\right)
\right\vert \neq1$. If this condition holds then $M_{i,j}$ is the \emph{only}
non-zero entry in the $i$-th row and the $j$-th column of $M\left(  D\right)
$. Then the form of $M\left(  D\right)  $ is
\[
\left[
\begin{array}
[c]{ccc}%
M^{\prime} &  & \\
& 1 & \\
&  & M^{\prime\prime}%
\end{array}
\right]  ,
\]
where the matrix $M^{\prime}$ is $\left(  i-1\right)  \times i$ and the matrix
$M^{\prime\prime}$ is $\left(  n-i+1\right)  \times\left(  n-i+2\right)  $. If
$D\in\mathcal{U}$ then $M^{\prime}$ and $M^{\prime\prime}$ have to be square.
Then $D\notin\mathcal{U}$. A contradiction.

Suppose that $\left\{  v_{i},v_{j}\right\}  $ is a bridge. Then $M_{i+1,i}=1$.
A similar reasoning as in the case of directed bridges applies. This forces
$M\left(  D\right)  $ to take one of the two forms
\[%
\begin{tabular}
[c]{lll}%
$\left[
\begin{array}
[c]{ccc}%
M^{\prime} &  & \\
&
\begin{array}
[c]{cc}%
0 & 1\\
1 & 0
\end{array}
& \\
&  & M^{\prime\prime}%
\end{array}
\right]  $ & and & $\left[
\begin{array}
[c]{ccc}%
M^{\prime} &  & \\
&
\begin{array}
[c]{cc}%
1 & 1\\
1 & 1
\end{array}
& \\
&  & M^{\prime\prime}%
\end{array}
\right]  $%
\end{tabular}
,
\]
where $M^{\prime}$ is $\left(  i-1\right)  \times\left(  i-1\right)  $ and
$M^{\prime\prime}$ is $\left(  n-1+2\right)  \times\left(  n-i+1\right)  $.
Clearly, $D\in\mathcal{U}$ if and only if $M^{\prime}$ and $M^{\prime\prime}$
support unitary matrices.

Suppose that $v_{i}$ is a cut-vertex. We can then label $V\left(  D\right)  $
such that $M\left(  D\right)  $ has the form
\[
\left[
\begin{array}
[c]{cccccc}
&  &  & 0 &  & \\
& M_{1} &  & \vdots & \mathbf{0} & \\
&  &  & 0 &  & \\
M_{i,1} & \cdots & M_{i,i} & M_{i,i+1} & \cdots & M_{i,n}\\
&  & M_{i+1,i} &  &  & \\
& \mathbf{0} & \vdots &  & M_{2} & \\
&  & M_{n,i} &  &  &
\end{array}
\right]  ,
\]
where $M_{i,1},...,M_{i,n},M_{i+1,i},...,M_{n,i}\in\left\{  0,1\right\}  $,
but not all are zero. Suppose that $D\in\mathcal{U}$. It is immediate to
observe that a similar reasoning as in the previous cases applies again,
forcing the cut-vertex to be in a connected component that is either
$\overleftrightarrow{K}_{2}$ or $\overleftrightarrow{K}_{2}^{+}$.
\end{proof}

\subsection{Corollaries}

Here we observe some corollaries of Theorem \ref{th1}.

\begin{corollary}
\label{co}Let $D$ be a connected graph on $n+2$ vertices and let
$D\in\mathcal{U}$. Then $D$ is $2$-vertex-connected and $2$-edge-connected.
\end{corollary}

\begin{corollary}
\label{indpath}Let $D$ be a connected graph on $n>2$ vertices and let
$D\in\mathcal{U}$. Then $D$ contains at least two independent paths between
any two vertices.
\end{corollary}

\begin{proof}
From the Global Version of Menger's theorem (see, \emph{e.g.}, \cite{D00},
Theorem 3.3.5).
\end{proof}

\begin{corollary}
Let $D$ be a connected graph on $n\geq3$ vertices and let $D\in\mathcal{U}$.
Then, for all $v_{i},v_{j}\in V\left(  D\right)  $, where $v_{i}\neq v_{j}$,
there exists a cycle containing both $v_{i}$ and $v_{j}$.
\end{corollary}

\begin{proof}
From Theorem 3.15 in \cite{M01}.
\end{proof}

A $k$\emph{-flow} in a graph $D$ is an assignment of an orientation of $D$
together with an integer $c\in\left\{  1,2,...,k-1\right\}  $ such that, for
each vertex $v_{i}$, the sum of the values of $c$ on the arcs into $v_{i}$
equals the sum of the values of $c$ on the arcs from $v_{i}$.

\begin{corollary}
\label{sey}Let $D$ be a graph. If $D\in\mathcal{U}$ then it has a $6$-flow.
\end{corollary}

\begin{proof}
From Seymour's theorem \cite{S81}.
\end{proof}

\begin{remark}
A \emph{cycle cover} of a graph $D$ is a set of cycles, such that every edge
of $D$ lies in at least one of the cycles. The length of a cycle cover is the
sum of the lengths of its cycles. If $D\in\mathcal{U}$ the it has a cycle
cover with length at most $\frac{\left|  A\left(  D\right)  \right|  }%
{2}+\frac{25}{24}\left(  V\left(  D\right)  -1\right)  $, since this result
holds for bridgeless graphs \cite{F97}.
\end{remark}

\begin{remark}
In a connected graph, a \emph{pendant-vertex} is a vertex with degree 1. The
graph of a unitary is without pendant-vertices. Zbigniew \cite{Z82} proved
that the probability that a random graph on $n$ vertices has no pendant
vertices goes to $1$ as $n$ goes to $\infty$. Is a random graph bridgeless?
\end{remark}

\begin{remark}
A graph $D$ is said to be \emph{eulerian} if $D$ is connected and the degree
of every vertex is even. An eulerian graph $D$ is said to be \emph{even} (odd)
if it has an even (odd) number of edges. Delorme and Poljak \cite{DP93} stated
the following conjecture which Steger confirmed for $d=3$: for $d\geq3$, every
bridgeless $d$-regular graph $D$ admits a collection of even eulerian
subgraphs such that every edge of $D$ belongs to the same number of subgraphs
from the collection. It would be interesting to verify if the conjecture is
confirmed by the graphs of unitary matrices.
\end{remark}

\subsection{Matchings}

The \emph{term rank} of a matrix is the maximum number of nonzero entries of
the matrix, such that no two of them are in the same row or column. Let
$M_{1}\circ M_{2}$ be the \emph{Hadamard product} of matrices $M_{1}$ and
$M_{2}$: $\left(  M_{1}\circ M_{2}\right)  _{ij}=M_{1_{ij}}M_{2_{ij}}$.

\begin{proposition}
\label{perm}Let $D$ be a digraph and let $D\in\mathcal{U}$. Then there exists
a permutation matrix $P$, such that $M\left(  D\right)  \circ P=P$.
\end{proposition}

\begin{proof}
If a digraph on $n$ vertices $D\in\mathcal{U}$ then the term rank of $M\left(
D\right)  $ is $n$. In fact, it is well-known that the possible maximum rank
of a matrix with digraph $D$ is equal to its term rank, that is the term rank
of $M\left(  D\right)  $. The proposition follows.
\end{proof}

A \emph{cycle factor} of a digraph $D$ is a collection of pairwise
vertex-disjoint dicycles spanning $\ D$. In other words, a cycle factor is a
spanning $1$-regular subdigraph of $D$.

\begin{proposition}
The digraph of a unitary matrix has at least a cycle factor
\end{proposition}

\begin{proof}
By Proposition \ref{perm}, since the adjacency matrix of a cycle factor is a
permutation matrix.
\end{proof}

\begin{remark}
The existence of a cycle factor and strong connectdness are necessary and
sufficient conditions for some families of digraphs to be hamiltonian (see,
\emph{e.g.}, \cite{B-JG01}).
\end{remark}

In a graph $D$, a \emph{perfect }$2$\emph{-matching} is a spanning subgraph
consisting of vertex-disjoint edges and cycle. A perfect $2$-matching is what
Tutte calls $Q$\emph{-factor} \cite{T53}.

\begin{proposition}
\label{2-m}Let $D$ be a graph without loops and let $D\in\mathcal{U}$. Then
$D$ has a perfect $2$-matching.
\end{proposition}

\begin{proof}
By Proposition \ref{perm}, there is a permutation matrix $P$ such that
$M\left(  D\right)  \circ P=P$. Since $M\left(  D\right)  $ is symmetric,
there is $P^{-1}$ such that $M\left(  D\right)  \circ P+P^{-1}=P+P^{-1}$.
Clearly, $P $ can be symmetric itself and in such a case $P\circ P^{-1}=P$.
The proposition follows. We consider a graph without loops, because in such a
case $P$ might have a fixed-point, and the fact $M\left(  D\right)  \circ P=P$
would not necessarily implies a perfect 2-matching.
\end{proof}

\begin{proposition}
\label{math}Let $D$ be a graph and let $D\in\mathcal{U}$. Then, for every
$S\subset V\left(  D\right)  $, $\left|  S\right|  \leq\left|  N\left(
S\right)  \right|  $.
\end{proposition}

\begin{proof}
Let $U$ be a unitary matrix acting on an complex vector space $\mathcal{H}$.
Since $U$ is invertible and since $U^{-1}=U^{\dagger}$, $U$ is an isomorphism
from $\mathcal{H}$ onto $\mathcal{H}$. The proposition follows, as a
consequence of the fact that an isomorphism is a bijective map.
\end{proof}

In a graph $D$ on $n=2k$ vertices, a \emph{matching} is collection of pairwise
vertex-disjoint graphs $\overleftrightarrow{K}_{2_{i}}$. If a matching has
$n/2$ members it is then called \emph{perfect matching}.

\begin{proposition}
Let $D$ be a bipartite graph and let $D\in\mathcal{U}$. Then $D$ has at least
a perfect matching.
\end{proposition}

\begin{proof}
By Proposition \ref{math}, together with the K\"{o}nig-Hall matching theorem
(see, \emph{e.g.}, \cite{LP86}).
\end{proof}

\begin{remark}
A graph $D$ has a perfect matching if and only if there exists a symmetric
permutation matrix $P$ without fixed points, such that $M\left(  D\right)
\circ P=P$. The conditions for the existence of perfect matchings in
non-bipartite graphs of unitary matrices remain to be studied.
\end{remark}

\subsection{A remark:\ perfect 2-matchings and the Sperner capacity of a
graph}

We observe now a consequence of Proposition \ref{2-m}. Consider a probability
measure $\mu$ with domain $V\left(  D\right)  $. Let $\left\{  v_{i}%
,v_{j}\right\}  $ be an edge of $D$. The vertices $v_{i}$ and $v_{j}$ induces
a subgraph $\overleftrightarrow{K}_{2}$. All the subgraphs of $D$ induced by
two connected vertices form the \emph{edge family} of $D$, denoted by
$\mathcal{F}\left(  D\right)  $. The \emph{entropy} of $\overleftrightarrow
{K}_{2}$ (see, \emph{e.g.}, \cite{GKV94}) is defined by
\[
H\left(  \left\{  v_{i},v_{j}\right\}  ,\mu\right)  =\left[  \mu\left(
v_{i}\right)  +\mu\left(  v_{j}\right)  \right]  \cdot h\left(  \frac
{\mu\left(  v_{i}\right)  }{\mu\left(  v_{i}\right)  -\mu\left(  v_{j}\right)
}\right)  ,
\]
where $h$ denotes the binary entropy function%
\[
h\left(  x\right)  =-x\log_{2}x-\left(  1-x\right)  \log_{2}\left(
1-x\right)  .
\]
The \emph{Sperner capacity} of $\mathcal{F}\left(  D\right)  $ \cite{CKS88}
is
\[
\Theta\left(  \mathcal{F}\left(  D\right)  \right)  =\max_{\Pr}\min_{\left\{
v_{i},v_{j}\right\}  \text{ in }D}H\left(  \left\{  v_{i},v_{j}\right\}
,\mu\right)  .
\]
This quantity has an information theoretical interpretation (it is related to
the zero-error capacity of channels) and it is used in the asymptotic solution
of various problems in extremal set theory (determination of the asymptotic of
the largest size of qualitative independent partitions in the sense of
R\'{e}nyi)\ \cite{GKV94}.

\begin{proposition}
Let $D$ be a graph and let $D\in\mathcal{U}$. Then $\Theta\left(
\mathcal{F}\left(  D\right)  \right)  =\frac{2}{n}$ and the corresponding
probability distribution is uniform over $V\left(  D\right)  $.
\end{proposition}

\begin{proof}
By Proposition \ref{2-m} and by Theorem 1 in \cite{G98}.
\end{proof}

\subsection{A conjecture about hamiltonian cycles}

Let $D\in\mathcal{U}$ be a connected graph on $n$ vertices. It is licit to ask
if the fact that $D\in\mathcal{U}$ is a sufficient condition for the existence
of hamiltonian cycles. Take as hypothesis the quadrangularity condition and
the existence of a perfect 2-matching. If $n=2,...,6$, it can be shown that
these two facts, together, imply the existence of an hamiltonian cycle.

\begin{conjecture}
Let $D$ be a connected graph and let $D\in\mathcal{U}$. Then $D$ is hamiltonian.
\end{conjecture}

\begin{remark}
A \emph{claw} is the bipartite graph $K_{1,3}$. Let \underline{$\lambda$} be
the graph with adjacency matrix
\[
\left[
\begin{array}
[c]{cccc}%
0 & 1 & 0 & 0\\
1 & 0 & 1 & 1\\
0 & 1 & 0 & 1\\
0 & 1 & 1 & 0
\end{array}
\right]  .
\]
We know that if $D$ is $2$-vertex-connected and if no induced subgraph of $D$
is isomorphic to $K_{1,3}$ or to \underline{$\lambda$}, then $D$ is
hamiltonian (see, \emph{e.g.}, \cite{M01}, Theorem 5.15). These conditions are
not sufficient to show that if $D\in\mathcal{U}$ then $D$ is hamiltonian. A
counterexample is the graph $D$ with adjacency matrix
\[
M\left(  D\right)  =\left[
\begin{array}
[c]{cccccc}%
0 & 0 & 0 & 1 & 1 & 0\\
0 & 0 & 0 & 1 & 1 & 1\\
0 & 0 & 0 & 1 & 1 & 1\\
1 & 1 & 1 & 0 & 0 & 0\\
1 & 1 & 1 & 0 & 0 & 0\\
0 & 1 & 1 & 0 & 0 & 0
\end{array}
\right]  .
\]
Since the matrix
\[
\left[
\begin{array}
[c]{ccc}%
\frac{1}{\sqrt{2}} & -\frac{1}{\sqrt{2}} & 0\\
\frac{1}{2} & \frac{1}{2} & \frac{1}{\sqrt{2}}\\
\frac{1}{2} & \frac{1}{2} & -\frac{1}{\sqrt{2}}%
\end{array}
\right]
\]
is unitary, $D\in\mathcal{U}$. It is easy to see that $D$ is hamiltonian even
if $D$ has a claw. Its adjacency matrix is
\[
\left[
\begin{array}
[c]{cccc}%
0 & 1 & 1 & 1\\
1 & 0 & 0 & 0\\
1 & 0 & 0 & 0\\
1 & 0 & 0 & 0
\end{array}
\right]  ,
\]
a submatrix of $M\left(  D\right)  $.
\end{remark}

\textbf{Acknowledgments} I would like to thank Scott Aaronson, Richard Jozsa
and Andreas Winter for reading and commenting earlier versions of this paper.
I\ would like to thank Peter Cameron for introducing me to Sonia Mansilla, and
Sonia for referring me to \cite{MS01}. This paper has been written while I was
a Student Member at the Mathematical Science Research Institute (MSRI),
Berkeley, CA, USA. I\ would like to thank MSRI for financial support.

\end{document}